\documentclass{crm-article}

\usepackage{citehack}
\usepackage{algorithm}
\usepackage{algorithmic}
\floatname{algorithm}{Алгоритм}

\newcommand{\argmin}{\mathop{\mathrm{argmin}}}
\newcommand{\Argmin}{\mathop{\mathrm{Argmin}}}
\newcommand{\LMO}{\mathop{\mathrm{LMO}}}
\newcommand{\Conv}{\mathop{\mathrm{Conv}}}

\newcommand*{\e}{\varepsilon}    
\newcommand{\sign}{\mathop{\mathrm{sign}}}

\begin{document}

\journalVol{10}

\journalNo{1}
\setcounter{page}{1}

\journalSection{Математические основы и численные методы моделирования}
\journalSectionEn{Mathematical modeling and numerical simulation}

\journalReceived{01.06.2016.}

\journalAccepted{01.06.2016.}


\affiliationnoref

\emailnoref

\UDC{519.85}
\title{Нижние оценки для методов типа условного градиента для задач минимизации гладких сильно выпуклых функций}
\titleeng{Lower bounds for conditional gradient type methods for minimizing smooth strongly convex functions
}

\author{\firstname{А.\,Д.}~\surname{Агафонов}}
\authorfull{Артём Дмитриевич Агафонов}
\authoreng{\firstname{A.\,D.}~\surname{Agafonov}}
\authorfulleng{Artem D. Agafonov}
\email{agafonov.ad@phystech.edu}
\affiliation{Национальный исследовательский университет «Московский физико-технический институт»,\protect\\ Россия, 141701, Московская облаcть,
г. Долгопрудный, 
Институтский пер., 9.}
\affiliationeng{National Research University Moscow Institute of Physics and Technology,\protect\\ 9 Institutskiy per., Dolgoprudny, Moscow Region, 141701, Russia}



\begin{abstract}
В данной работе рассматриваются методы условного градиента для оптимизации сильно выпуклых функций. Это методы, использующие линейный минимизационный оракул, то есть умеющие вычислять решение задачи 
$\Argmin_{x\in X}{\langle p,x \rangle}$
для заданного вектора $p \in \mathbb{R}^n$. Существует целый ряд  методов условного градиента, имеющих линейную скорость сходимости в сильно выпуклом случае. Однако, во всех этих методах в оценку скорости сходимости входит размерность задачи, которая в современных приложениях может быть очень большой. В данной работе доказывается, что в сильно выпуклом случае скорость сходимости методов условного градиента в лучшем случае зависит от размерности задачи $n$ как $\widetilde{\Omega}(\sqrt{n})$. Таким образом, методы условного градиента могут оказаться неэффективными для решения сильно выпуклых оптимизационных задач больших размерностей.

Отдельно рассматривается приложение методов условного градиента к задачам минимизации квадратичной формы.
Уже была доказана эффективность метода Франка-Вульфа для решения задачи квадратичной оптимизации в выпуклом случае на симплексе (PageRank). Данная работа показывает, что использование методов условного градиента  для 
минимизации квадратичной формы в сильно выпуклом случае малоэффективно из-за наличия размерности в оценке скорости сходимости этих методов. Поэтому рассматривается метод рестартов условного градиента (Shrinking Conditional Gradient). Его отличие от методов условного градиента заключается в том, что в нем используется модифицированный линейный минимизационный оракул, который  для заданного вектора $p \in \mathbb{R}^n$ вычисляет решение задачи 
$$\Argmin\{\langle p, x \rangle: x\in X, \|x-x_0\| \leqslant R \}.$$ 
В оценку скорости сходимости такого алгоритма размерность уже не выходит. С помощью рестартов метода условного градиента получена сложность (число арифметических операций) минимизации квадратичной формы на $\infty$-шаре. Полученная оценка работы метода сравнима со сложностью градиентного метода.

\end{abstract}

\begin{abstracteng}
In this paper, we consider conditional gradient methods for optimizing strongly convex functions. These are methods that use a linear minimization oracle, which, for a given vector $p \in \mathbb{R}^n$, computes the solution of the subproblem 
$$\Argmin_{x\in X}{\langle p,x \rangle}.$$
There are a variety of conditional gradient methods that have a linear convergence rate in a strongly convex case. However, in all these methods, the dimension of the problem is included in the rate of convergence, which in modern applications can be very large. In this paper, we prove that in the strongly convex case, the convergence rate of the conditional gradient methods in the best case depends on the dimension of the problem $ n $ as $ \widetilde {\Omega} (\sqrt {n}) $.  Thus, the conditional gradient methods may turn out to be ineffective for solving strongly convex optimization problems of large dimensions.

Also, the application of conditional gradient methods to minimization problems of a quadratic form 
is considered. 
The effectiveness of the Frank-Wolfe method for solving the quadratic optimization problem in the convex case on a simplex (PageRank) has already been proved. This work shows that the use of conditional gradient methods to solve 
the minimization problem of a quadratic form in a strongly convex case is ineffective due to the presence of dimension in the convergence rate of these methods. Therefore, the Shrinking Conditional Gradient method is considered. Its difference from the conditional gradient methods is that it uses a modified linear minimization oracle. It's an oracle, which, for a given vector $p \in \mathbb{R}^n$, computes the solution of the subproblem 
$$\Argmin\{\langle p, x \rangle: x\in X, \|x-x_0\| \leqslant R \}.$$ 
The convergence rate of such an algorithm does not depend on dimension. Using the Shrinking Conditional Gradient method  the complexity (the total number of arithmetic operations) of solving the minimization problem of quadratic form on a $ \infty $-ball is obtained. The resulting evaluation of the method is comparable to the complexity of the gradient method.

\end{abstracteng} 
\keywordeng{Conditional Gradient method}
\keywordeng{Frank-Wolfe method}
\keywordeng{Shrinking Conditional Gradient}

\maketitle

\paragraph{Введение}
В данной работе рассматриваются методы оптимизации, использующие линейный минимизационный оракул. Это такой оракул, который для данного вектора $p \in \mathbb{R}^n$ вычисляет решение задачи
$$\Argmin_{x\in X}{\langle p,x \rangle}.$$
Такими методами являются классический метод условного градиента (Франка-Вульфа)  \cite{frank1956wolfe} и различные его варианты \cite{jaggi2013revising}. Эти методы -- одни из самых популярных методов условной оптимизации вследствие их простоты и отсутствия проектирования на шаге алгоритма  \cite{dvurechensky2020selfconcordant}. В основе метода условного градиента лежит идея линеаризации целевой функции в точке и минимизация линейной аппроксимации на бюджетном множестве. Однако, этот метод сходится достаточно медленно, сублинейно ($O(1/\e)$) \cite{levitin1966}. 

Ряд недавних работ показал, что в случае сильно выпуклых функций оценку можно улучшить до линейной ($O(\log 1/\e)$) \cite{garber2013}, \cite{kerdreux2019}, \cite{lacoste2015global_conv}. Однако, во всех алгоритмах, представленных в работах выше, в оценку скорости сходимости входит размерность задачи (см. дискуссию в \cite{lacoste2015global_conv}). Стоит отметить, что такая проблема возникает только в сильно выпуклом случае. Например, в оценку скорости сходимости классического метода Франка-Вульфа, который работает для выпуклых задач, размерность не входит. 

В данной работе мы показываем, что для задачи сильно выпуклой оптимизации в оценку скорости сходимости методов, использующих линейный минимизационный оракул, обязательно входит размерность задачи как $\tilde{\Omega}(\sqrt{n})$. В современных приложениях размерность может быть очень большой, что негативно сказывается на времени работы метода. Таким образом, методы условного градиента могут оказаться неэффективными для решения задач больших размерностей в сильно выпуклом случае.

Также рассматривается приложение методов условного градиента к задаче квадратичной оптимизации на единичном шаре в $\infty$-норме. Такая постановка задачи на $\infty$-шаре, или, что эквивалентно, параллелепипеде встречается достаточно часто \cite{gornov2009book}.  Например, она появляется в задачах, где, известно, что каждая компонента ограниченна, то есть принадлежит заданному отрезку. 
В данной работе рассматривается случай, когда матрица $A$ дважды разреженна (одновременно разреженна по столбцам и строкам)  \cite{anikin2015sparse}. Если считать, что матрица $A$ имеет размер $n \times n$ , а число элементов в каждой строке и столбце не больше чем $s \ll n $, то число ненулевых элементов в матрице может быть $sn$.

Было показано \cite{anikin2015sparse}, \cite{gasnikov2016diss}, что классический метод Франка-Вульфа эффективно решает задачу
минимизации квадратичной формы  на симплексе (PageRank \cite{brin1998page}).
В данной работе рассматривается возможность  решения задачи квадратичной оптимизации методами типа условного градиента в сильно выпуклом случае. Алгоритмы, использующие линейный минимизационный оракул, не могут решать эту задачу эффективно. В данной работе доказывается, что оценки скорости сходимости таких методов содержат размерность. Для решения задачи минимизации квадратичной формы авторами используется метод рестартов условного градиента (Shrinking Conditional Gradient) \cite{lan2013lo}. Вместо линейного минимизационного оракула используется оракул, вычисляющий решение задачи
$$\Argmin\{\langle p, x \rangle: x\in X, \|x-x_0\| \leqslant R \}.$$ 
В сложность этого метода размерность уже не входит. Анализ работы этого метода показывает, что сложность метода (число арифметических операций):
$$\frac{8L}{\mu}O(n)+\frac{8L}{\mu}\left(sn\log\frac{\mu R_0}{\e}\right),$$
где первый член -- затраты на препроцессинг. Эта оценка сравнима с сложностью градиентного метода \cite{allen2014}. 

Статья состоит из введения и трех основных разделов. В первом разделе рассматривается классический метод условного градиента (Франка-Вульфа) и вводятся необходимые вспомогательные понятия. Второй раздел посвящен вопросам размерности в оценке скорости сходимости методов условного градиента. В третьем разделе предлагается подход к решению задачи минимизации квадратичной формы на $\infty$-шаре с помощью рестартов метода условного градиента. 

\paragraph{Классический метод условного градиента}




Рассмотрим следующую задачу оптимизации:
\begin{equation}\label{opt_task}
    f(x) \rightarrow \min \limits_{x \in X},
\end{equation}
где 
\begin{itemize}
    \item $X \in \mathbb{R}^n$ -- выпуклое и замкнутое множество,
    \item $f(x)$ -- выпуклая функция,
    \item $f(x)$ ограничена снизу на $X$ и достигает своего минимума в $x^*$.
\end{itemize}

Рассмотрим классический метод условного градиента \cite{frank1956wolfe}. В его основе лежит идея линеаризации функции. На $k$-ом шаге метода в точке $x_k$  линеаризуем функцию $f(x)$, минимизируем эту линейную аппроксимацию на $X$ и найденную точку используем как направление движения.

\begin{algorithm}
\caption{Класический метод условного градиента (Франка-Вульфа)}
\begin{algorithmic}[1]\label{alg:classicCNdG}
\STATE \textbf{Дано:} $x_0 \in X$ -- начальная точка. Положим $y_0 = x_0$.
\FOR{$k = 1, 2, \ldots$} 
\STATE Вычисляем $x_k = \argmin_{x\in X}{\langle f'(y_{k-1}),x \rangle}.$
\STATE Положим $y_k = y_{k-1} + (1-\gamma_k)x_k$, где $\gamma_k \in [0,1]$.
\ENDFOR
\end{algorithmic}
\end{algorithm}

Чтобы гарантировать сходимость классического метода условного градиента необходимо правильно выбрать длину шага  $\gamma_k$. Рассмотрим две стратегии выбора $\gamma_k$:
\begin{enumerate}
    \item Последовательность $\{\gamma_k\}$ выбирается заранее:
    \begin{equation}\label{step1}
    \gamma_k = \frac{2}{k+1}, \text{ где } k = 1, 2, \ldots.
    \end{equation}
    \item Наискорейший спуск:
    \begin{equation}\label{step2}
    \gamma_k = \argmin\limits_{\gamma \in [0,1]} f(y_{k-1} + (1-\gamma)x_k).
    \end{equation}
\end{enumerate}

Выделим класс функций с липшицевым градиентом. 
\begin{fed}
    Будем говорить, что функция $f$ имеет непрерывный по Липшицу (или просто липшицев) градиент, если  
    \begin{equation}
        \|\nabla f(x) - \nabla f(y)\|_* \leqslant L \|x - y\|,~\forall x, y \in X.
    \end{equation}
\end{fed}
Пусть у нас есть выпуклая функция с липшицевым градиентом. Известно, что классическому методу условного градиента с длиной шага, выбранной по правилу \eqref{step1} или \eqref{step2}, требуется 
\begin{equation}\label{classic_cmplxty}
    O\left(\dfrac{LD^2_X}{\e}\right)
\end{equation} 
итераций, чтобы достигнуть точности $\e$ по функции, где $D_x = \max \limits_{x, y \in X}\|x-y\|$  \cite{frank1956wolfe}. Заметим, что несмотря на то, что в методе условного градиента нам не нужно выбирать норму, в оценку сложности метода входят зависящие от нее константы $L = L_{\|\cdot\|}, D_X = D_{X, \|\cdot\|}$. Так как оценка \eqref{classic_cmplxty} верна для произвольной нормы, то сложность данного метода можно оценить точнее:
\begin{equation}\label{classic_cmplxty2}
    O\left(\inf \limits_{\|\cdot\|}\dfrac{L_{\|\cdot\|}D^2_{X, \|\cdot\|}}{\e}\right).
\end{equation}

\paragraph{Размерность в методах условного градиента}



Будем рассматривать различные методы условного градиента. Классический метод Франка-Вульфа является одним из них. Все эти методы итеративно используют линейный минимизационной оракул для решения задачи \eqref{opt_task}. 
\begin{fed}\label{def:LMO}
    Для данного вектора $p \in \mathbb{R}^n$ линейный минимизационный оракул ($\LMO(p)$) вычисляет решение следующей задачи:
    \begin{equation}\label{eq:LMO}
        \LMO(p) \in \Argmin_{x\in X}{\langle p,x \rangle}.
    \end{equation}
\end{fed}

В общем виде алгоритм методов условного градиента можно записать следующим образом  \cite{lan2013lo}.

\begin{algorithm}
\caption{Общий вид методов условного градиента}  
\begin{algorithmic}[1]\label{alg:genLCP}
\STATE \textbf{Дано:} $x_0 \in X$ -- начальная точка. 
\FOR{$k = 1, 2, \ldots$} 
\STATE Определяем вектор $p_k \in \mathbb{R}^n$.
\STATE Вычисляем $x_k = \LMO(p_k)$.
\STATE Возвращаем $y_k = \Conv\{x_0, \ldots, x_k\}$.
\ENDFOR
\end{algorithmic}
\end{algorithm}

Заметим, что в Алгоритме \ref{alg:genLCP} не указаны конкретное правило вычисления $y_k$ и стратегия выбора $p_k$. Это делает алгоритм достаточно общим. Например, Алгоритм \ref{alg:genLCP} включает в себя, как частный случай, классический метод условного градиента (Алгоритм \ref{alg:classicCNdG}).

Рассмотрим задачу \eqref{opt_task} с функцией $f$ с липшицевым градиентом. Для методов условного градиента известна следующая нижняя оценка\cite{lan2013lo}. Число итераций Алгоритма \ref{alg:genLCP} для решения этой задачи с точностью $\e$ по функции не превосходит 
\begin{equation}\label{eq:low_bnd}
    \left\lceil \min \left\{\dfrac{n}{2}, \dfrac{LD^2}{4\e}\right\}\right\rceil - 1.
\end{equation}

Пусть в сильно-выпуклом случае сложность некоторого метода условного градиента $\mathcal{M}$ (Алгоритм \ref{alg:genLCP}) есть $O\left(\max\left(\ln^p\frac1\e, 1\right)\right)$, где $p = 1, 2, \ldots$. Например, при $p = 1$ скорость сходимости метода является линейной, а при $p >= 2$ -- сублинейной. Ускорим метод условного градиента с помощью  Каталиста(Catalyst) \cite{lin2018}, \cite{gasnikov2018book}.


Покажем, как получается ускорение. Введем функции 
\begin{gather*}
    F_{\varkappa, x} (y) = f(y) + \frac{\varkappa}{2}\|y - x\|_2^2,\\
    f_\varkappa(x) = \min \limits_{y \in X} F_{\varkappa, x}(y) = F_{\varkappa, x}(y_\varkappa(x)).
\end{gather*}
Меняем исходную задачу \eqref{opt_task} на следующую:
\begin{equation}\label{eq:catalyst}
    f_{\varkappa} (x) \rightarrow \min \limits_{x \in X},
\end{equation}
где $\varkappa$ мы можем выбирать сами.

Конструкция Каталист состоит из двух вложенных циклов. Рассмотрим метод Монтейро-Свайтера, использующийся в Каталисте для ускорения методов \cite{gasnikov2018book}, \cite{monteiro2013}, \cite{ivanova2019}.  В данном случае удобно выбрать $\varkappa = L$.

\begin{algorithm}
\caption{Метод Монтейро-Свайтера}  
\begin{algorithmic}[1]\label{alg:catalyst}
\STATE \textbf{Дано:} 

\STATE \textbf{Инициализация:} $z_0, y_0, A_0 = 0$
\WHILE{критерий остановки не выполнен} 
\STATE Вычисляем 
    $$a_{k+1} = \frac{1/\varkappa + \sqrt{1/\varkappa^2 + 4A_k/\varkappa}}{2}, ~ A_{k+1} = A_k + a_{k+1}$$

\STATE Вычисляем 
$$x_{k+1} = \frac{A_k}{A_{k+1}}y_k+\frac{a_{k+1}}{A_{k+1}}z_k$$
\STATE Подбираем $y_{k+1}$ так, чтобы выполнялось условие Монтейро-Свайтера
\begin{equation}\label{eq:monteiro}
    \|\nabla F_{\varkappa, x_{k+1}} (y_{k+1}) \|_2 \leqslant \frac{\varkappa}{2}\|y_{k+1} - x_{k+1}\|_2
\end{equation}
\STATE Вычисляем $$z_{k+1} = z_{k} - a_{k+1}\nabla f(y_{k+1}) $$
\ENDWHILE
\end{algorithmic}
\end{algorithm}

Условие Монтейро–Свайтера \eqref{eq:monteiro} позволяет вместо точного решения $y_\varkappa (x)$ воспомогательной задачи \eqref{eq:catalyst} искать неточное решение задачи \cite{gasnikov2018book}, \cite{lin2018}
\begin{equation}\label{eq:catalapprox}
    y_{k+1} \approx \argmin F_{\varkappa, x_{k+1}} (x).
\end{equation}
Таким образом, на каждом шаге внутреннего цикла методом $\mathcal{M}$ решается задача \eqref{eq:catalapprox} чтобы найти $y_{k+1}$, удовлетворяющий \eqref{eq:monteiro}. Сложность решения этой задачи c точностью $\e$ равна $O\left(\max\left(\ln^p\frac1\e, 1\right)\right).$  

Во внешней итерации конструкции Каталист используется Алгоритм \ref{alg:catalyst}. Его сложность равна $O\left(\sqrt{\frac{LR^2}{\e}}\right)$ \cite{monteiro2013}, \cite{gasnikov2018book}. Тогда итоговое время работы метода окажется следующим  $$O\left(\sqrt{\frac{LR^2}{\e}}\max\left(\ln^p\frac1\e, 1\right)\right).$$
Сравним его с нижней оценкой \eqref{eq:low_bnd}

$$\sqrt{\frac{LR^2}{\e}}\max\left(\ln\frac1\e, 1\right) \geqslant  \min \left\{ \frac n2, \frac{LR^2}{4\e}  \right\}.$$
Предположим, что  в данный момент  метод работает в режиме  $ \frac{LR^2}{2\e} \approx n $.
Пусть $C_n$ -- константа в оценке сложности метода, то есть метод сходится как $C_n\ln^p\frac1\e = O(\ln^p\frac1\e)$. Получаем 
$$C_n \geqslant \dfrac{\sqrt{\frac{LR^2}{\e}}}{\ln^p\frac{LR^2}{\e}} = \frac{\sqrt{n}}{\ln^p{n}} = \widetilde{\Omega}(\sqrt n).$$
Таким образом, мы доказали следующую  теорему.

\begin{teo}
Пусть у  оптимизационного метода $\mathcal{M}$, имеющего структуру Алгоритма \ref{alg:genLCP}, сложность в сильно-выпуклом случае равна $C_n\left(\max\left(\ln^p\frac1\e, 1\right)\right)$, где $C_n$ -- некоторая константа, $p = 1, 2, \ldots$, а $C_n$ -- некоторая константа. Тогда $C_n = \widetilde{\Omega}(\sqrt{n})$.
\end{teo}

Получается, что в оценке скорости сходимости методов типа условного градиента в сильно выпуклом случае нельзя избавиться от размерности.

\paragraph{Методы условного градиента для решения задачи квадратичной оптимизации на кубе}







Из результатов предыдущего пункта следует, что методы условного градиента (Алгоритм \ref{alg:genLCP}) будут неэффективны для решения  многих сильно выпуклых оптимизационных  задач, так как в оценке скорости сходимости таких алгоритмов нельзя избавиться от размерности. Ланом был предложен метод рестартов условного градиента (Shrinking Conditional Gradient)\cite{lan2013lo} для сильно выпуклых функций. В нем вместо линейного минимизационного оракула (Определение \ref{def:LMO}) используется линейный оракул, способный решать оптимизационные задачи следующего вида
\begin{equation}\label{eq:ELMOmain}
    \min\{\langle p, x \rangle: x\in X, \|x-x_0\| \leqslant R \}.
\end{equation}

\begin{algorithm}
\caption{Рестарты условного градиента (Shrinking Conditional Gradient)}
\begin{algorithmic}[1]\label{alg:ShCndG}
\STATE \textbf{Дано:} $x_0 \in X$ -- начальная точка, $R_0 = D_X$. 
\FOR{$t = 1, 2, \ldots$} 
\STATE Положим $y_0 = p_{t-1}$
\FOR{$k = 1, \ldots, \frac{8L}{\mu}$}
\STATE Вычисляем $x_k$ -- решение задачи:
\begin{equation}\label{eq:ELMO}
 \Argmin_{x\in X_{t-1}} \langle f'(y_{k-1}),x\rangle, \text{ где }  X_{t-1} = \{x \in X: \|x - p_{t-1}\| \leqslant R_{t-1}\}.
 \end{equation}
\STATE Положим $y_k = (1 - \gamma_k)y_{k-1} + \gamma_k x_k$ для некоторого $\gamma_k \in [0, 1]$.
\ENDFOR
\STATE Положим $p_t = y_k$ и $R_t = \frac{R_{t-1}}{\sqrt{2}}$
\ENDFOR
\end{algorithmic}
\end{algorithm}

Стратегия выбора длины шага $\gamma_k$ в Алгоритме \ref{alg:ShCndG} такая же, как и в Алгоритме \ref{alg:classicCNdG}, и задается формулами \eqref{step1}, \eqref{step2}. Рестартам условного градиента с такими правилами выбора шага требуется не более чем 
\begin{equation} 
    \frac{8L}{\mu}\left\lceil \max \left( \frac{\mu R^2_0}{\e}, 1\right) \right\rceil
\end{equation}
итераций, чтобы достичь точности $\e$ по функции \cite{lan2013lo}.


Рассмотрим задачу квадратичной оптимизации на кубе ($\infty$-шаре):
\begin{equation}\label{eq:PR}
    f(x) = \frac12 \|Ax\|^2_2 \rightarrow \min \limits_{\|x\|_\infty \leqslant 1},
\end{equation}
 При этом мы считаем, что в каждом столбце и каждой строке матрицы $A$ не более $s \ll n$ элементов отлично от нуля ($A$ -- разрежена), то есть всего в $A$ может быть $sn$ ненулевых элементов. Кроме того, $A$ удовлетворяет условию $\mu I \preccurlyeq A$, то есть задача \eqref{eq:PR} сильно выпукла.

Для решения этой задачи будем пользоваться Алгоритмом \ref{alg:ShCndG}. Сначала рассмотрим шаги внутреннего цикла этого метода (строки $5$-$6$ Алгоритма \ref{alg:ShCndG}). Пусть в \eqref{eq:ELMO} $X_{t-1} = \{x: \|x\|_\infty \leqslant 1\}$. Решаем задачу  
\begin{equation}\label{eq:PRtask}
    \Argmin \limits_{\|x\|_\infty \leqslant 1} \langle f'(y_{k}),x\rangle.
\end{equation} 
Так как это задача нахождения минимума линейной формы на кубе, то ее решение --  вершина этого куба. Нам известен градиент с предыдущего шага $\nabla f(y_{k-1}) = (\partial f(y_{k-1})/ \partial x^1, \ldots, \partial f(y_{k-1})/ \partial x^n)^T$. Тогда решением задачи \eqref{eq:PRtask} является вектор
\begin{equation*}
    x_k = - \left(\sign\left(\frac{\partial f(y_{k-1})}{\partial x^1}\right), \ldots, \sign\left(\frac{\partial f(y_{k-1})}{\partial x^n}\right)\right)^T.
\end{equation*}
Сложность нахождения этого вектора есть $O(n)$.

Пересчитаем $y$ 
\begin{equation*}
    y_k = (1 - \gamma_k)y_{k-1} + \gamma_k x_k,
\end{equation*} 
и градиент:
\begin{equation*}
    \nabla f(y_k) = A^TAy_k = (1 - \gamma_k)A^TA y_{k-1} + \gamma_k A^TAx_k.
\end{equation*} 
Сложность пересчета градиента -- $O(s^2n)$, так как сначала надо умножить разряженную матрицу $A$ на $n$-мерный вектор $x_k$, а потом еще домножить на $A^T$.

Пусть мы вышли из внутреннего цикла Алгоритма \ref{alg:ShCndG}. Чтобы решить задачу \eqref{eq:ELMO} необходимо найти $X_t$ -- пересечение исходного множества $X$ и шара в $\infty$-норме с центром в $p_t$ и радиусом $R_t$.  Опишем процедуру его получения.

На каждой итерации внешнего цикла Алгоритма \ref{alg:ShCndG} мы будем делать препроцессинг: масштабировать систему координат и переносить её центр, чтобы в пересечении получался шар $B_1^\infty(0)$. Так всегда можно сделать, потому что мы пересекаем $n$-мерные осепараллельные прямоугольники. 

Введем обозначение:
\begin{equation*}
    B_R^\infty (c) \stackrel{\text{def}}{=} \{x: \|x - c\|_\infty \leqslant R\}.
\end{equation*}

После масштабирования и изменения центров системы координат на предыдущих итерациях, $X$ имеет вид осепараллельного прямоугольника. Пусть его центр находится в точке $c_t$, и задан вектор расстояний от центра до граней $a_t$ (или координаты одной вершины, по которым можно вычислить эти расстояния). Проверим, лежит ли $B^\infty_{R_t} (p_t)$ в $X$. Это можно сделать за $O(n)$
, проверив принадлежность отрезка $[p_t^i - R_t,  p_t^i + R_t]$ отрезку $[c_t^i - a_t^i,  c_t^i + a_t^i]$ для всех $i \in \{1, 2, \ldots, n\}$.
 
Пусть оказалось, что $B^\infty_{R_t} (p_t) \subseteq X$. Перенесем $p_t$ в нуль. Обновим значение $y_0 = p_{t} - p_{t} = 0$. Сделаем масштабирование системы координат, чтобы привести шар $B^\infty_{R_t} (p_t) \subseteq X$ к единичному. 
Осталось пересчитать  множество $X$ в соответсвии с произведенными преобразованиями:
$$c_{t+1} = \frac{c_t - p_t}{R_t}, ~~~ a_{t+1} = \frac{a_t}{R}.$$

В другом случае, $B^\infty_{R_t} (p_t) \nsubseteq X$. В этом случае для всех $i \in [1, 2, \ldots, n]$ ищем пересечение $[p_t^i - R_t,  p_t^i + R_t]  \cap [c_t^i - a_t^i,  c_t^i + a_t^i] = [\min^i, \max^i]$ и центры $m_t^i$ этих пересечений.  Вектор $m_t$ -- центр осепараллельного параллелепипеда, образованного пересечением $B^\infty_{R_t} (p_t)$ и $X$. Перенесем $m_t$ в ноль  и  отмасштабируем по каждой из осей, чтобы сделать $X_{t}$ единичным кубом. Следовательно, координаты любой точки $x \in \mathbb{R}^n$ преобразуются следующим образом:  
\begin{equation}
    x^i := \frac{2 \cdot (x^i - m_t^i)}{\max^i - \min^i}.
\end{equation}
Начальная точка внутреннего цикла $y_0$ и центр $c_t$ исходного множества $X$  преобразуются в соответсвии с этим правилом.
Сложность этих вычислений есть $O(n)$.

Таким образом, мы свели задачу \eqref{eq:ELMO} к \eqref{eq:PRtask}. Процедуру решения задачи \eqref{eq:PRtask} мы описали выше.

В итоге сложность шагов внутреннего цикла Алгоритма $\ref{alg:genLCP}$ равна $O(sn)$, сложность препроцессинга -- $O(n)$. Получаем следующую оценку сложности предложенного метода:
\begin{equation}
    \frac{8L}{\mu}O(n)+\frac{8L}{\mu}\left(s^2n\log\frac{\mu R_0}{\e}\right),
\end{equation}
где $R_0 = 2$ -- размер шара $B_1^\infty(0)$ в $\infty$-норме, а $L$ -- константа Липшица градиента в $\infty$-норме.

Сравним эту оценку с  оценкой градиентного метода \cite{allen2014}:
\begin{equation}
    \sqrt{\frac{L}{\mu}} O\left(s^2n\log\frac{L_1R_0^2}{\e}\right),
\end{equation}
где $R_0^2 = O(\sqrt{n})$, если выбрать прокс-функцию $d(x) = \frac{1}{2}\|x\|_2^2$ \cite{gasnikov2018book}. 

Получается, что предложенный нами метод решения задачи \eqref{eq:PR} работает не сильно лучше градиентного метода. Заметим, что не получилось улучшения оценки, как в статье \cite{anikin2015sparse}, где методом Франка-Вульфа (Алгоритм \ref{alg:classicCNdG}) решали задачу поиска вектора RageRank в выпуклом случае на единичном симплексе.  Дело в том, что если задача поставлена на симплексе, то сложность пересчета градиента, при хранении его компонент в бинарной куче, занимает $O(s^2\log n)$, так как необходимо обновлять всего одну компоненту градиента на каждой итерации. В случае $\infty$-шара надо обновлять все $n$ компонент. Соответсвенно, сложность этой операции выше ($O(s^2n)$). 

Удобство постановки задачи на $\infty$-шаре заключается в том, что на каждой итерации сохраняется геометрия задачи. Пересечение $\infty$-шаров, на котором идет минимизация в задаче \eqref{eq:ELMO}, есть многомерный осепараллельный прямоугольник. Стоит отметить, что пока неясно  как использовать Алгоритм \ref{alg:ShCndG} для решения задачи PageRank на симплексе или $1$-шаре. Пересечение $1$-шара и шара в произвольной норме, или симплекса и шара, может оказаться сложным выпуклым многогранником. Непонятно, как найти этот многогранник. То есть,  уметь , например, вычислять все его вершины (в какой-то из них и будет минимум линейного функционала в задаче \eqref{eq:ELMO}). Несмотря на то, что авторы статьи \cite{lan2013lo}
указывают в ней, что сложность решения задачи \eqref{eq:ELMOmain} сравнима со сложностью решения задачи \eqref{eq:LMO}, пример про задачу PageRank на симплексе или $1$-шаре показывает противное.    

С другой стороны, в Алгоритме \ref{alg:ShCndG}  в задаче \eqref{eq:ELMO} не обязательно брать пересечение исходного множества с шаром, а можно брать некоторое множество из исходного, которое содержит в себе пересечение. Это повлияет только на константу в оценке скорости сходимости метода. Таким образом, на каждой итерации метода можно всегда решать задачу \eqref{eq:ELMO} на одном и том же множестве, например, $1$-шаре. Но и в этом случае непонятно, как именно найти множество, содержащее $X_{t}$.


\dotfill

Автор выражает благодарность Гасникову Александру Владимировичу за полезные обсуждения и идеи.

\end{document}